\documentclass[12pt]{article}
\usepackage[latin1]{inputenc}
\usepackage{amssymb}
\newtheorem{theorem}{Theorem}
\newtheorem{proposition}[theorem]{Proposition}
\newtheorem{lemma}[theorem]{Lemma}

\newtheorem{corollary}[theorem]{Corollary}

\newtheorem{remark}[theorem]{Remark}

\newtheorem{example}[theorem]{Example}

\newcommand{\R}{\mathbb{R}}

\newcommand{\spa}{\mbox{span}}
\newcommand{\Ric}{\mbox{Ric}}

\newcommand{\rank}{\mbox{rank}}
\newcommand{\kerl}{\mbox{ker}}

\newcommand{\po}{{\hspace*{-1ex}}{\bf .  }}
\newcommand{\ii}{isometric immersion }

\def\P{{\cal P}}
\def\Sal{{\cal S}}
\def\<{\langle}

\def\>{\rangle}
\def\a{\alpha}

\newcommand{\N}{\mathcal{N}}
\newcommand{\nap}{\nabla^\perp}
\newcommand{\nab}{\tilde\nabla}
\def\bea{\begin{eqnarray*} }
\def\eea{\end{eqnarray*} }
\def\be{\begin{equation} }
\def\ee{\end{equation} }
\def\nap{\nabla^\perp}
\def\proof{\noindent{\it Proof: }}
\def\qed{\ifhmode\unskip\nobreak\fi\ifmmode\ifinner
\else\hskip5 pt \fi\fi\hbox{\hskip5 pt \vrule width4 pt
height6 pt  depth1.5 pt \hskip 1pt }}
\begin{document}

\title{Submanifolds with nonparallel first normal\\ bundle revisited 
\thanks {{\it Mathematics Subject Classification 2000 
 53B25.}}}
\author{ Marcos Dajczer and Ruy Tojeiro\footnote{
Partially supported by  FAPESP grant
11/21362-2.}}
\date{}
\maketitle

\begin{abstract} 
In this paper, we analyze the  geometric structure of a Euclidean  submanifold whose
osculating spaces form a nonconstant family of proper subspaces 
of the same  dimension.  We prove that if the rate of change of the osculating spaces  
is small, then  the submanifold must be a (submanifold of a)  ruled submanifold of a 
very special type. We also give a sharp estimate of the dimension of the rulings.

\end{abstract}

The osculating space of a Euclidean submanifold $M^n$ at a point is the
subspace of Euclidean space that is spanned by the tangent and curvature vectors of
all smooth curves in $M^n$ through that point. 
If all osculating spaces along $M^n$ coincide with a fixed subspace $H$, it is 
an elementary fact that $M^n$ is contained in an affine subspace parallel to
$H$. Thus, it is a natural problem to study for which submanifolds the osculating
spaces form a nonconstant family of proper subspaces of the same  dimension.
 In  this paper, we show that if the rate of change of the osculating 
spaces  is small, in a sense to be made precise below, 
then the submanifold must be contained in a ruled submanifold of a very special type.

Let $f\colon  M^n\to \R^N$ denote an \ii of an $n$-dimensional connected 
Riemannian manifold into Euclidean space.  The \emph{first normal space} of $f$
at $x\in M^n$ is the normal subspace $N_1^f(x)\subset N_fM(x)$  spanned by its
second fundamental form  $\a_f$, that is,
$$
N_1^f(x)=\spa\{\a_f(X,Y) : X,Y\in T_xM\}.
$$
The \emph{osculating space} of $f$ at $x\in M^n$ is defined as $f_*T_x M\oplus N_1^f(x)$. 
It is easy to see that all  osculating spaces of $f$ have the same dimension and are parallel  to a 
fixed proper subspace of $\R^N$ if and only if the  first normal spaces form a proper normal 
subbundle $N_1^f$ that is parallel in the normal connection; see \cite{Da} or  
\cite{Sp}. Then $f$ reduces codimension to $p=\rank \,N_1^f$, 
that is, it can be seen as a substantial \ii into an affine subspace $\R^{n+p}$ of
$\R^N$.

A rather simple argument  shows that  $N_1^f$ must be parallel in the normal
connection if $p<n$ and at any $x\in M^n$ the $s$-nullities $\nu_s$ of $f$ satisfy  
\be\label{condition}
\nu_s(x)<n-s
\ee
for all $1\leq s\leq p$; see \cite{Da}, \cite{Dr} or (\ref{d1}) below. 
Recall that
$$
\nu_s(x)=\max_{U^s\subset N_1^f(x)}\dim\N(\a_{U^s})
$$
where $U^s\subset N_1^f(x)$ is any $s$-dimensional vector subspace and
$$
\N(\a_U(x))=\{Y\in T_xM: \a_U(Y,X)=0\;\mbox{for all}\; X\in T_xM\}
$$
for $\a_U=\pi_U\circ\a_f$ and $\pi_U\colon N_1^f\to\,U$  the orthogonal projection.
Notice that $\nu_p(x)$ is the standard \emph{index of relative nullity}
$\nu_f(x)=\dim\N(\a_f(x))$, that is, the dimension of the \emph{relative nullity}
subspace of $f$ at $x\in M^n$.

Consider the subspace $\Sal(x)$ of $N_1^f(x)$ spanned by the projections 
onto $N_1^f(x)$ of the derivatives $\nab_X\mu$ in the ambient space,  with 
$X\in T_xM$, of local sections $\mu\in (N_1^f)^\perp$ of its orthogonal complement 
in the normal bundle $N_fM$. If all subspaces $\Sal(x)$ have the same  dimension  
along $M^n$, and thus form a vector subbundle $\Sal=\Sal_f$,  we
may say that  the rank $s$ of $\Sal$ measures to what extent the first 
normal bundle $N_1^f$ fails to be parallel.

If $\Sal$ coincides with $N_1^f$ and $p\leq 6$, it turns out that condition 
(\ref{condition}) fails for the relative nullity, i.e.,  $\nu_f\geq n-p>0$ at any point. 
The latter has strong well-known geometric consequences, namely, the submanifold  carries 
a $\nu_f$-dimensional totally  geodesic foliation whose leaves are open subsets of 
affine subspaces in $\R^N$.

Our main result is that there is a single class of submanifolds for which  $\Sal$  is 
a proper subbundle of $N_1^f$ of rank $s\leq 6$, any other example being 
a submanifold of an element of this class. These are ruled submanifolds, with rulings
of dimension at least $n-s$, for which $\Sal$ is constant in the ambient space along 
the rulings. In particular, the rulings  belong to the kernel of $\a_\Sal$, 
and therefore condition (\ref{condition}) is violated for $s$. Examples of such submanifolds, 
showing  that the preceding estimate on the dimension of the rulings is sharp, 
are constructed in the last section.

As discussed in the next section,  the results of this paper   
generalize  those in \cite{Dt} for $p\leq 3$. We also point out  that, although   
stated for submanifolds of Euclidean space, our results can easily be extended to
ambient spaces of constant sectional curvature.

\section{The result}

  In this section, we first give a precise statement of our main result and 
then discuss some particular cases.\vspace{1.5ex}

Let $f\colon M^n\to\R^N$ denote a locally substantial isometric immersion of a 
connected  Riemannian manifold, i.e., there is no open subset $U\subset M^n$ such that 
$f(U)$ is  contained in a proper affine subspace of $\R^N$. Assume that $f$ is 1-regular, i.e., 
the first normal spaces $N_1^f(x)$ have constant dimension $p$.  Thus, these subspaces form 
a vector subbundle $N_1^f$ of the normal bundle $N_fM$ which we assume to be proper, i.e., 
$p<N-n$.

Assume $p<n$ and let $\phi\colon (N_1^f)^\perp\oplus TM\to N_1^f$ be the tensor
defined  by
$$
\phi(\mu,X)=(\nap_X\mu)_{N_1^f}
$$
where $(\;)_{N_1^f}$ denotes the $N_1^f$-component.
We say that $f$ has \emph{nonparallel first normal bundle at $x\in M^n$} if   
$\phi(x)\neq 0$, i.e., if the dimension $s(x)$ of the normal  vector subspace 
$\Sal(x)\subset N_1^f(x)$ given by 
$$
\Sal(x)=\spa\{\phi(\mu,X):\mu\in (N_1^f)^\perp(x)\;\mbox{and}\;X\in T_xM\}
$$
is nonzero. Thus, along each connected component of the open dense subset 
of $M^n$ where  $s(x)=s$ is constant, the vector subspaces $\Sal(x)$
form a   vector subbundle $\Sal$ of $N^f_1$.\vspace{1ex}

In the following statement, that an \ii $F\colon  N^m\to\R^N$, $m>n$, is an 
\emph{extension} of the isometric immersion $f\colon M^n\to\R^N$  means that 
there exists an isometric embedding $i\colon  M^n \to N^m$ such that  $f=F\circ i$. 
Also, by $f$ being $d$-ruled we understand that there exists a $d$-dimensional 
integrable distribution in $M^n$ whose leaves are (mapped by $f$ into) open 
subsets of affine subspaces in the ambient space.

\begin{theorem}\po\label{nonparallel}
Let $f\colon M^n\to\R^N$ be a 1-regular locally substantial isometric immersion such
that $s(x)=s$ is constant with  $0<s<n$ and $s\leq 6$. Then, either
\begin{itemize}
\item[(i)] $s=p$ and $f$ has index of relative nullity $\nu_f\geq n-p$,
 or 
\item[(ii)] $1=s<p$ and $f$  has an extension $F\colon  N^{n+p-1}\to\R^N$
such that  $\nu_F= n+p-2$ and $N^F_1$ is nonparallel of rank one, or
\item[(iii)] $1<s<p$ and there is an open dense subset of $M^n$, the union of 
open subsets $U_{k,d}$ with  $d\geq n-s$ and $n-d\leq k\leq q:=n-d+p-s$, such that:
\item[(a)] $f|_{U_{q,d}}$  is $d$-ruled and $\Sal_f$  is constant in $\R^N$  along the 
rulings, and
\item[(b)] $f|_{U_{k,d}}$, $k<q$,  has a ruled extension $F\colon N^{n+q-k}\to\R^N$ 
such that $N^F_1$ is nonparallel of rank $p+k-q$ and $\Sal_F$ is constant along the rulings. 
The rulings have dimension $n+p-k-s$ and coincide with $\N(\a_F)$  if $k=n-d$.  

Moreover, if $s=2$ we have that  $U_{k,d}=\emptyset$ for $k\geq 5$. 
\end{itemize}
\end{theorem}

The ruled extensions  in parts $(ii)$ and $(b)$ of $(iii)$ are as in $(i)$ 
and $(a)$ of $(iii)$, respectively.  

  For a ruled Euclidean submanifold, it is easily seen that for any vector $X$ tangent 
to a ruling  the Ricci curvature  satisfies $\Ric(X)\leq 0$, with equality if and 
only if $X$ belongs to the relative nullity subspace.  Hence, we have
the following immediate consequence of  Theorem \ref{nonparallel}.  
  
\begin{corollary}\po Under the assumptions of Theorem \ref{nonparallel}, cases $(i)$  
and $(iii)\!-\!(a)$ cannot occur if $\Ric_M>0$. If  $\Ric_M\geq 0$ then $f|_{U_{q,d}}$ 
in case $(iii)\!-\!(a)$ satisfies $\nu_f=d$.
\end{corollary}

To illustrate Theorem \ref{nonparallel}  we discuss next the cases $p=1,2$ and $3$. 
Notice that these are the cases that have already been considered in \cite{Dr}. 

\begin{example}\label{first}\po{\em The case $p=1$. Here, the only possibility 
is that $s=1$, and hence $\nu_f = n-1$. In particular, the manifold  $M^n$ is flat.
}\end{example}

Submanifolds as  above can be easily described parametrically.  For instance, consider the 
image under the normal exponential map of a parallel normal subbundle of the normal bundle of 
a curve with non-vanishing curvature; see also Theorem 1 in \cite{Dt}.

\begin{example}\label{second}\po{\em The case $p=2$.   We only have the following
two  possibilities:
\begin{itemize}
\item[(i)]  $s=2$, and hence  $\nu_f = n-2$. 
\item[(ii)]  $s=1$, in which case  $f$ admits an extension  
$F\colon N^{n+1}\to\R^N$ such that  $\nu_F = n$ (hence $N^{n+1}$ is flat) 
and $N^F_1$ is nonparallel of rank one. 
\end{itemize}
}\end{example}

The submanifolds in case $(i)$ have been studied in \cite{Df} and \cite{Dm}, 
where a  parametric classification has been obtained in most cases.

\begin{example}\label{third}\po{\em The case $p=3$.  
Then one of the following holds: 
\begin{itemize}
\item[(i)]   $s=3$ and $f$ satisfies $\nu_f\geq n-3$. 
\item[(ii)]  $s=1$ and  $f$ has an extension $F\colon N^{n+2}\to\R^N$ 
such that $\nu_F = n+1$ ($N^{n+2}$ is flat) and $N^F_1$ is  nonparallel of rank one. 
\item[(iii)]   $s=2<k=3$, in which case $f$ is $(n-2)$-ruled and $\Sal$ is constant
along the rulings.
\item[(iv)]   $s=2=k$ and $f$ has an extension $F\colon N^{n+1}\to\R^N$ 
such that $\nu_F = n-1$ and $N^F_1$ has rank two.
\end{itemize}
}\end{example}

Observe that $F$  in $(ii)$ of Example \ref{second} and Example \ref{third} is as 
$f$ in Example~\ref{first}. Also, the extension $F$  in  $(iv)$  of Example \ref{third} is 
as $f$ in $(i)$ of Example \ref{second}.

\section{A class of ruled extensions}

In this section of independent interest, we find sufficient conditions for an
Euclidean submanifold to admit a ruled  extension carrying a normal subbundle
that is constant in the ambient space along the rulings. We point out that
a special case was already considered in \cite{Dt}.
\vspace{1.5ex}

Let $f\colon  M^n\to \R^N$ be an \ii satisfying the following conditions:
\begin{itemize}
\item[(i)] Its normal bundle splits  
orthogonally and
smoothly into two  vector subbundles
$$
N_fM=L\oplus P
$$
such that the rank $\ell$ of $L$ satisfies $0<\ell< N-n$. 
\item[(ii)]  The   
subspaces
$$
D(x)=\N(\a_P(x))\subset T_xM
$$
have constant dimension $d>0$ on $M^n$ (thus form a tangent  
subbundle $D\subset TM$). 
\item[(iii)] The subbundle $P$ is parallel along $D$  in the  normal connection, 
thus in $\R^N$. Hence, also $L$ is parallel along $D$  in the  normal connection.         
\end{itemize}

Let $\gamma\colon E\oplus P\to E\oplus L$ be the tensor given by
\be\label{gamma}
\gamma(Y,\mu)=(\nab_Y\mu)_{E\oplus L}=-A_\mu Y + (\nap_Y\mu)_L,
\ee
where the subbundle $E\subset TM$ of rank $n-d$ is defined by the orthogonal splitting 
$TM=D\oplus E$ and $\nab$ denotes the connection in $\R^N$.                                             

At $x\in M^n$, let $\Gamma(x)\subset E(x)\oplus L(x)$ be the subspace defined by
\be\label{eq:gamma}
\Gamma(x)=\spa\{\gamma(Y,\mu): Y\in E\;\mbox{and}\;\mu\in P\}.
\ee
Since $E$ is spanned by the vectors $A_\mu Y$ for $\mu\in P$ and $Y\in  
E$, it follows
from (\ref{gamma}) that
\be\label{cases}
n-d\leq \dim \Gamma(x)\leq n-d+\ell.
\ee
Assume further that
\begin{itemize}
\item[(iv)]
$\dim\Gamma(x)=k$ is constant on $M^n$. 
\end{itemize}
 Let  
$\pi\colon \Lambda\to M^n$
be the affine vector bundle of rank $r = n-d+\ell-k$ that is defined by the  
orthogonal splitting
$$
\Gamma^k\oplus\Lambda^r=E^{n-d}\oplus L^\ell.
$$

\begin{lemma}\po\label{par}
The distribution  $D$ is integrable and $\Lambda\cap TM=\{0\}$ holds.
\end{lemma}

\proof Take $\mu\in P$ and $Z,Y\in D$. Since $P$ is parallel along   
$D$ in $\R^N$, we have from
\be\label{curtensor}
0=\widetilde{R}(Y,Z)\mu=\nab_Y\nab_Z\mu-\nab_Z\nab_Y\mu-\nab_{[Y,Z]}\mu
\ee
that  $\widetilde\nabla_{[Y,Z]}\mu\in P$. Hence $A_\mu[Y,Z]=0$, and  
thus $D$ is integrable.

Take $Z\in \Lambda\cap TM$. Then $Z\in E$ and
$$
0=\<Z,\nab_X\mu\>=-\<A_\mu Z,X\>
$$
for any $\mu\in P$ and $X\in TM$. Thus $Z\in D$ and hence  
$Z=0$.\qed\vspace{1,5ex}

The affine subspaces $\Delta (x)$ defined by
$$
\Delta(x)=D(x)\oplus\Lambda(x)
$$
form an affine bundle over $M^n$ of rank $d+r = n+\ell-k$.

\begin{lemma}\label{rulings}\po
The bundle  $\Delta$ is parallel in $\R^N$  along the leaves of $D$.
\end{lemma}

\proof  It suffices to show that the orthogonal complement  
$\Gamma\oplus P$ of $\Delta$ in $\R^N$                               
is parallel in $\R^N$ along the leaves of $D$.  First observe that      
$$
\Gamma\oplus P =\spa\{\nab_X\mu: X\in TM\;\mbox{and}\;\mu\in P\}.
$$
Then, we have from (\ref{curtensor}) that
$$
\nab_Y\nab_X\mu=\nab_X\nab_Y\mu+\nab_{[Y,X]}\mu\in \Gamma\oplus P
$$
for any $\mu\in P$, $Y\in D$ and $X\in TM$, and the assertion follows.  
\qed\vspace{1,5ex}

Define $F\colon  N^{n+r}\to\R^N$ as the restriction of the map
$$
\lambda\in\Lambda\mapsto f(\pi(\lambda))+ \lambda
$$
to a tubular neighborhood  $N^{n+r}$ of the $0$-section
$j\colon M^n\hookrightarrow N^{n+r}$ of $\Lambda$ where it is an immersion.
Then $f=F\circ j$ and
\be\label{tangent}
T_{j(x)}N=j_*T_xM\oplus\Lambda(x)
\ee
for any $x\in M^n$.

Lemma \ref{rulings} yields that $F$ is ruled with
$\Delta(\lambda):=\Delta(\pi(\lambda))$ as the ruling through $\lambda\in \Lambda$. 
For $\lambda\in\Lambda$, $\mu\in P$ and $X\in TM$, it follows from                      
$$
\<\nab_X\lambda,\mu\>=- \<\lambda,\nab_X\mu\>=0
$$
that $\P\subset N_F N$ where $\P(\lambda)=P(\pi(\lambda))$. Moreover,  we have that
$$
\Delta=\N(\a^F_{\P}).
$$
In fact, the inclusion $\Delta\subset\N(\a^F_{\P})$ holds because $\P$  
is constant along $\Delta$.
For the opposite inclusion observe that $\a_\P^F|_{{TM\times TM}}=\a_{P}$.
We easily obtain from (\ref{tangent}) that equality
is satisfied along $M^n$. To conclude the proof observe that the  
dimension of $\N(\a^F_{\P})$
can only decrease along $\Delta\subset N^{n+r}$ from its value on
$M^n$ if $N^{n+r}$ is taken small enough.
\vspace{1,5ex}

We summarize the above facts in the following statement.

\begin{proposition}\label{FF}\po Let $f\colon M^n\to \R^N$ be an isometric 
immersion satisfying  $(i)$-$(iv)$ above. Then $f$ admits a ruled extension  
$F\colon N^{n+r}\to\R^N$, $r = n-d+\ell-k$, with the following properties:
\begin{itemize}
\item[(a)] The distribution $\Delta$ of rulings of $F$ satisfies  
$D^d(x)=\Delta^{d+r}(x)\cap T_xM$ at any $x\in M^n$.
\item[(b)] There is an orthogonal splitting $N_FN=\mathcal{L}\oplus\P$ so that 
$\rank\;\mathcal{L}=\ell-r$, $\Delta=\N(\a^F_{\P})$ and $\P$ is constant in $\R^N$ 
along $\Delta$.
\end{itemize}
Moreover, we have: 
\begin{itemize}
\item [(c)] If $r=0$ then $f$ is $d$-ruled and $P$ is constant 
in $\R^N$ along the rulings.         
\item [(d)] If $r=\ell$ then $\Delta$ is the relative nullity distribution of $F$.              
\end{itemize}
\end{proposition}

\section{The proof}

A key ingredient in  the proof of Theorem \ref{nonparallel} is a basic property of 
regular elements of a bilinear form observed by Moore \cite{Mo}. It is stated below 
as Proposition \ref{moore}.  \vspace{1.5ex}

Let $\beta\colon  V\times U\to W$ be a  bilinear form  between finite
dimensional real vector spaces. We call $Z\in V$   a
(left) {\it regular element} of $\beta$ if the map
$\beta_Z={\beta}(Z,\,\cdot\,)\colon  U\to W$
satisfies
$$
\dim \beta_Z(U)=\max\{\dim \beta_Y(U):\, Y\in V\},
$$
and denote by $RE(\beta)$ the subset of regular elements of $\beta$.
It is a well-known fact that the set $RE(\beta)$ is open and dense in $V$.

\begin{proposition}\po\label{moore}
If $\beta\colon  V\times U\to W$ is a  bilinear form and $Z\in RE(\beta)$, then
$$
\beta(V,\ker\beta_Z)\subset \beta_Z(U).
$$
\end{proposition}

With the notations from Section $1$,   consider a 1-regular locally substantial 
isometric immersion $f\colon M^n\to\R^N$  such that $s(x)$ has a constant 
value $0<s< n$.

\begin{lemma}\po \label{D1} It holds that  $\N(\phi)=\N(\a_\Sal)$.
\end{lemma}

\proof Let $\mu_1\in RE(\phi)$ be a globally defined unit vector field and set
$\phi_{\mu_1}=\phi(\mu_1,\,\cdot\,)$.
Without loss of generality, we may assume that the subspaces $S_1(x)\subset \Sal(x)$ defined by
$$
S_1(x)=\phi_{\mu_1}(T_xM)
$$  
have constant dimension $1\leq s_1\leq s$.  Hence, the tangent subspaces
$D_1(x)=\kerl\phi_{\mu_1}(x)$ satisfy $\dim D_1(x) = n-s_1$. It suffices to show that
\be\label{d1}
D_1=\N(\a_{S_1}),
\ee
i.e.,  that $Y\in D_1$ if and only if $A_{\nap_X\mu_1}Y=0$ for any $X\in TM$. 
But this follows from the Codazzi equation 
\be\label{codazzi}
A_{\nap_X\delta}Y=A_{\nap_Y\delta}X
\ee
for any $\delta\in (N_1^f)^\perp$.\qed

\begin{lemma}\po\label{Dnontrivial} Suppose that $s\leq 6$. Then
$D=\N(\phi)$  satisfies
\be\label{claim}
\dim D\geq n-s.
\ee
\end{lemma}
\proof
Let $\mu_1$ be as in the previous lemma. Again, we may assume that 
$S_1(x)$ has constant dimension $1\leq s_1\leq s$ on $M^n$.
In view of Lemma \ref{D1}, the assertion holds if $s_1=s$. If $s_1< s$,
consider the orthogonal splitting
$$
\Sal=S_1\oplus S_1^\perp
$$
and let $\psi\colon (N_1^f)^\perp\oplus TM\to S_1^\perp$ denote the bilinear form
defined by
$$
\psi(\mu,X)=(\nap_X\mu)_{S_1^\perp}.
$$
Take  $\mu_2\in RE(\phi)\cap RE(\psi)$ and set $t=\dim\psi(\mu_2,TM)$.
Then $S_2=\phi_{\mu_2}(TM)$ satisfies
$$
\dim\, (S_1+S_2)=s_1+t\;\;\mbox{and}\;\;\dim  S_1\cap S_2=s_1-t.
$$
It follows using Proposition \ref{moore} that
\be\label{4}
\dim D_1\cap D_2\geq \dim D_1-\dim S_1\cap S_2\geq n-2s_1+t.
\ee
 If $t=s_1$ then $S_1\cap S_2=0$. Thus $D_1=D_2$.  
In particular (\ref{claim}) holds if $s_1=1$ since this forces $t=1$. Therefore, 
we may assume
\be\label{1}
s_1\geq 2.
\ee

We first analyze the case $t=1$. In this case, we have that 
$H=\kerl\,\psi({\mu_2}, \cdot)$ is a hyperplane in
$TM$.  From (\ref{codazzi}) we obtain
$$
A_{\nap_Z\mu_2}X=A_{\nap_X\mu_2}Z=0
$$
for any $Z\in D_1$ and $X\in H$.  This implies that $\dim\phi_{\mu_2}(D_1)\leq
1$.
Otherwise, there would exist a two-dimensional plane in $S_1$ such that the corresponding  
shape operators would have the same kernel of codimension one.
But then a vector in this plane would belong to $(N_1^f)^\perp$, and this is a
contradiction.
It follows that $\dim D_1\cap D_2\geq n-s_1-1$.

If $\Sal=S_1+S_2$ then (\ref{claim}) holds
since $s=s_1+1$  and $D=D_1\cap D_2$. If otherwise, we just repeat the process
and obtain
subspaces $S_1, \ldots, S_m$ and $D_1, \ldots, D_m$, $m=s-s_1+1$, such that
$\Sal=S_1+\cdots +S_m$ and  $\dim D_1\cap\cdots \cap D_m \geq  
n-s_1-m+1-s$. Then
$D=D_1\cap\cdots \cap D_m$, and (\ref{claim}) follows. 

By the above, we may assume that $t\geq 2$.
We argue for the case $s=6$, the other cases being similar and easier.  If
$t=s_1$ then $s_1=2,3$. In these cases we have seen that $D_1=D_2$, and thus
(\ref{claim}) holds. Hence, in view of (\ref{1}) and $t\geq 2$ we may assume
that $s_1>t\geq 2$.
Thus, it remains to consider the cases $(s_1,t)=(3,2)$ and $(s_1,t)=(4,2)$. In
the latter case, we have that $\Sal=S_1+S_2$, and (\ref{claim}) follows from
(\ref{4}). In the first case, we have
$\dim\,(S_1+S_2)=5$, $\dim  S_1\cap S_2=1$ and $\dim D_1\cap D_2\geq n-4$. We
now repeat the process and obtain $S_3$ such that $\Sal=S_1+S_2+S_3$ and 
$\dim  S_i\cap S_j=1$ if $i\neq j$. In this case, it is now clear
that $\dim D\geq n-5$.\qed

\begin{remark}\po{\em Our proof does not work for $s=7$. In fact, in this case
we may have $s_1=5$ and $t=2$.  Thus $\Sal=S_1+S_2$ and  (\ref{4}) only
yields $\dim D\geq n-8$.
}\end{remark}

Now consider the global smooth orthogonal splitting $N^f_1=L^{p-s}\oplus \Sal^s$. 
Then, we have the global orthogonal splitting
\be\label{split}
N_fM=L^{p-s}\oplus P
\ee
where $P=\Sal^s\oplus (N_1^f)^\perp$.

\begin{lemma}\po \label{parallel} The subbundle  $P$ is parallel along $D$ 
in the normal connection.
\end{lemma}

\proof By the Ricci equation, we have 
$$
\nap_Y\nap_X\mu_1-\nap_X\nap_Y\mu_1-\nap_{[Y,X]}\mu_1=0.
$$
Take $Y\in D_1$ and $X\in TM$.  Then,
$$
\nap_Y(\nap_X\mu_1)_{S_1} + \nap_Y(\nap_X\mu_1)_{(N_1^f)\perp} 
=\nap_X\nap_Y\mu_1+\nap_{[Y,X]}\mu_1\in P.
$$
By Proposition \ref{moore}, the second term on the left-hand-side belongs to $P$. 
It follows that $\nap_Y\delta\in P$
for any $Y\in D_1$ and $\delta\in S_1$. 
\qed\vspace{1,5ex}

\noindent\emph{Proof of Theorem \ref{nonparallel}.}  Assume first that $s=p$, that is,  
that $\Sal=N_1^f$. Then, Lemma~\ref{D1} and Lemma~\ref{Dnontrivial}  imply that 
$\nu_f\geq n-p$. 

Suppose now that $s<p$. For each positive integer $d$, let $U_d$ denote the interior 
of the subset of all $x\in M^n$ such that the subspace  $D(x)$  has dimension $d$. 
It follows from Lemma~\ref{Dnontrivial} that $d\geq n-s$. 
By the lower semi-continuity of the dimension, we have that $\cup_{d} U_d$
is (open and) dense in $M^n$.
Now let $U_{k,d}$ be the interior of the subset of all $x\in U_d$ such that the 
subspace $\Gamma(x)$  given by   (\ref{eq:gamma}), with respect to the splitting (\ref{split}),
has dimension $k$. Then (\ref{cases}) with $\ell=p-s$ gives $n-d\leq k\leq q$. 
Again by the lower semi-continuity of the dimension, we have that $\cup_{ k} U_{k,d}$
is (open and) dense in $U_d$.
 
In view of Lemma \ref{D1} and Lemma \ref{parallel}, we can apply  
Proposition \ref{FF} for $f|_{U_{k,d}}$. If $k=q$, 
we obtain from Proposition \ref{FF} -$(c)$ that  $f|_{U_{q,d}}$    
is $d$-ruled and $P$ (hence $\Sal$)  is constant in $\R^N$  along the 
rulings. 

If $k<q$, it follows from Proposition \ref{FF} that $f$ admits a ruled extension 
$F\colon  N^{n+r}\to\R^N$, $r = n-d+\ell-k=q-k$,
with rulings of dimension $n+\ell-k = n+p-k-s$. Moreover,
there is an orthogonal splitting $N_FN=\mathcal{L}\oplus\P$,
where $\P$ is the parallel extension (in $\R^N$) of $P$ along the rulings, 
such that $\rank\;\mathcal{L}=p-s-r$. In particular, 
$\rank\;N_1^F=p-r=p+k-q$.

Finally, if  $k = n-d$ then  the rulings of $F$ 
coincide with  its relative nullity distribution by Proposition \ref{FF}-$(d)$.

The global assertion in $(ii)$ for the case $1=s<p$ is due to the fact that 
$s=1$ implies $d=1$, and also $k=1$, as follows from (\ref{gamma}). It is also a 
consequence of (\ref{gamma}) that $k\leq 4$ if $s=2$, hence in this case  
$U_{k, d}=\emptyset$ for $k\geq 5$.\qed

\section{Examples}

In this section we give examples of Euclidean  submanifolds  satisfying the 
conditions in part $(iii)-(a)$ of Theorem \ref{nonparallel}. More precisely, we 
construct ruled submanifolds $M^{2m}$ in $\R^{2m+6}$ with four dimensional first
normal bundle such that $\Sal$ has rank two and is constant along the codimensional 
two rulings. These examples show that the result cannot be improved since the 
rulings are not in the relative nullity distribution and their dimension
achieve the minimum possible value given by the estimate.\vspace{1,5ex}

Let $g\colon L^2\to\R^{2(m+3)}$, $m\geq 2$, be a substantial  elliptic surface
in the sense of \cite{Df}, i.e.,  there exists a (unique up to sign) almost
complex structure $J$ on $L^2$ such that
$$
\alpha_g(Z,Z)+\alpha_g(JZ, JZ)=0
$$
for any $Z\in TL$. 
For instance, the surface can be minimal, which is equivalent to  $J$ being orthogonal. 
Then, it turns out that the normal bundle  of $g$ splits orthogonally as
$$
N_gL=N_1^g\oplus\cdots\oplus N_{m+2}^g
$$
where each plane bundle $N_k^g$, $1\leq k\leq m+2$, is its $k^{th}$-normal bundle; 
see \cite{Df} for details. 
Recall that the \emph{$k^{th}$-normal space} $N_k^h$, $k\geq 2$, of an
isometric immersion $h\colon  M^n \to \R^N$ at $x\in M^n$  is defined as
$$
N_k^h(x)=\spa\{\alpha_h^{k+1}(X_1,\ldots, X_{k+1}) \,:\,X_1,\ldots, X_{k+1}\in T_xM\},
$$
where $\alpha_h^{\ell}\colon  TM\times \cdots \times TM\to N_hM$, $\ell \geq 3$, is the 
$\ell^{th}$-fundamental form given by 
$$
\alpha_h^{\ell}(X_1,\ldots, X_{\ell})
=\pi^{\ell-1}(\nabla^\perp_{X_\ell}\cdots \nabla^\perp_{X_3}\alpha(X_2,X_1)).
$$ 
Here  $\pi^\ell$ is the orthogonal projection onto 
$(N_1^h\oplus\ldots \oplus N_{\ell -1}^h)^\perp \cap N_hM$. 

Define $f\colon  M^{2m}\to\R^{2(m+3)}$ as the restriction of the map
$$
\xi\in N_1^g\oplus\cdots\oplus N_{m-1}^g\mapsto g(\pi(\xi))+ \xi
$$
to a tubular neighborhood of the $0$-section  $L^2$ of 
$\pi\colon N_1^g\oplus\cdots\oplus N_{m-1}^g\to L^2$
where it is an immersion. Given $\xi\in M^{2m}\setminus L^2$, we claim that
$$
f_*T_\xi M\oplus N_1^f(\xi)=g_*T_xL\oplus N_1^g(x)\oplus\cdots\oplus N_{m+1}^g(x),
\,\,\,\,x=\pi(\xi).
$$

Let $\tilde \xi$ be a local section of $N_1^g\oplus\cdots\oplus N_{m-1}^g$
on a neighborhood $U$ of $x$ such that $\tilde \xi(U)\subset M^{2m}$ and 
$\tilde\xi(x)=\xi$. Then 
\be\label{f*}
f_*\tilde \xi_*X=g_*X+\tilde \nabla_X\tilde \xi
\ee
for any $X\in T_xL$. On the other hand, for a vertical vector $V\in T_\xi M$ we have
$$
f_*V=V.
$$
Hence $N_1^g(x)\oplus\cdots\oplus N_{m-1}^g(x)\subset f_*T_\xi M$ and 
$f_*T_\xi M\subset g_*T_xL\oplus N_1^g(x)\oplus\cdots\oplus N_m^g(x)$.
Regarding the local section $\tilde \xi$ as a vertical vector field of $M^{2m}$, we obtain
\be\label{vec1}
\nab_X\tilde \xi=\nab_{\tilde \xi_*X}f_*\tilde \xi\in f_*T_\xi M\oplus N_1^f(\xi).
\ee
Thus $N_m^g(x)\subset f_*T_\xi M \oplus N_1^f(\xi)$, hence also  
$g_*T_xL\subset f_*T_\xi M \oplus N_1^f(\xi)$ by (\ref{f*}).
Differentiating (\ref{f*}) yields
$$
\nab_{\tilde \xi_*Y}f_*\tilde \xi_*\tilde X=\nab_Yg_*\tilde X+\nab_Y\nab_X\tilde\xi
$$
for all $X,Y\in T_xL$, where $\tilde X$ is any vector field  on a neighborhood  
of $x$ with $\tilde X(x)=X$.
Thus $N_1^f(\xi)\subset g_*T_xL\oplus N_m^g(x)\oplus N_{m+1}^g(x)$ and 
$N_{m+1}^g(x)\subset N_1^f(\xi)$,  
and the claim follows. 

Note also that the rulings of $f$ are not in its relative nullity distribution. In fact, 
it follows from (\ref{vec1}) that 
\be\label{n1}
\spa\{\alpha_f(Z,V): Z,V\in T_\xi M\,\mbox{and}\,\,V\,\,
\mbox{vertical}\}=(g_*T_xL\oplus N_{m}^g(x))\cap N_1^f(\xi).
\ee

 We have from the claim that $N_fM=N_1^f\oplus N_{m+2}^g$.
Thus, the immersion $f$  is ruled by  $N_1^g\oplus\cdots\oplus N_{m-1}^g$
and $\Sal=N_{m+1}^g$ has rank  two and is constant in the ambient 
space along the rulings. Moreover, by (\ref{n1}) the rulings are not in the relative nullity
distribution  and their dimension satisfy the equality in the estimate given
in part $(iii)-(a)$ of Theorem~\ref{nonparallel}.

{\renewcommand{\baselinestretch}{1}

\hspace*{-20ex}\begin{tabbing} \indent\= IMPA -- Estrada Dona Castorina, 110
\indent\indent\=  Universidade Federal de S\~{a}o Carlos\\
\> 22460-320 -- Rio de Janeiro -- Brazil  \>
13565-905 -- S\~{a}o Carlos -- Brazil \\
\> E-mail: marcos@impa.br \> E-mail: tojeiro@dm.ufscar.br
\end{tabbing}}
\end{document}